\input amstex.tex
\documentstyle{amsppt}
\magnification=1200

\loadbold

\def\Co#1{{\Cal O}_{#1}}
\def\Fi#1{\Phi_{|#1|}}
\def\fei#1{\phi_{#1}}
\def\rw{\rightarrow}
\def\lrw{\longrightarrow}

\def\Bbbp1{{\Bbb P}^1}
\def\Div{\text{\rm Div}}

\def\dim{\text{\rm dim}}
\def\roundup#1{\ulcorner{#1}\urcorner}
\def\rounddown#1{\llcorner{#1}\lrcorner}
\def\Bbbq{\Bbb Q}

\TagsOnRight

\topmatter
\title
Inequalities of Noether type for 3-folds of general type
\endtitle
\author
{Meng Chen}
\endauthor
\address \ \ \ \
\newline Institute of Mathematics, Fudan University
\newline Shanghai, 200433, PR China
\newline  {\it E-mail address:} {\bf meng\@x263.net}
\endaddress
\thanks
\endthanks
\thanks
\endthanks
\thanks
This paper is supported by the National Natural Science Foundation of China (No. 10131010),
Shanghai Scientific $\&$ Technical Commission (Grant 01QA14042) and SRF for ROCS, SEM.
\endthanks

\abstract
If $X$ is a smooth complex projective 3-fold with ample canonical divisor $K$, then the
inequality  $K^3\ge \frac{2}{3}(2p_g-7)$
holds, where $p_g$ denotes the geometric genus. This inequality is nearly sharp. We also give
similar, but more complicated, inequalities for general minimal 3-folds of general type.
\endabstract

\endtopmatter

\document
\vcorrection{-0.8cm}
\leftheadtext{M. Chen}
\rightheadtext{Inequalities of Noether type}

\head {\rm Introduction} \endhead
Given a minimal surface $S$ of general type, we have two famous inequalities, which play crucial
roles in detailed analysis of surfaces. One is the Bogomolov-Miyaoka-Yau inequality
$K_S^2\le 9\chi(S)$ (\cite{M1}, \cite{Y1}, \cite{Y2}), while the other is the classical Noether
inequality
$K_S^2\ge 2p_g-4\ge 2\chi(X)-6$. The fundamental importance of these inequalities in mind, M. Reid
asked in 1980s

\proclaim{Question 1} What would be the right analogue of the Noether inequality in dimension three?
\endproclaim

Let $X$ be a minimal threefold. If $K_X$ is Cartier and very ample, then
$K_X^3\ge 2p_g-6$ by Clifford's theorem applied to the intersection curve cut out by two general
members of  $|K_X|$. In 1992, Kobayashi \cite{Kob} studied Gorenstein canonical 3-folds and obtained
an effective, but partial, upper bound of $K_X^3$ in terms of $p_g(X)$ for such varieties. One of
his discoveries is that too naive a generalization of the classical Noether inequality is in general
false; there are a series of smooth projective 3-folds $X$ with ample canonical divisor such that
$$K_X^3=\frac{2}{3}(2p_g(X)-5),\ (p_g(X)=7,\ 10, \ 13,\cdots). \tag{0.1}$$
In what follows, we show that Kobayashi's examples indeed attain the minima of $K_X^3$, provided $X$
is smooth and $K_X$ is ample:

\proclaim{Corollary 2} If $X$ is a smooth complex projective 3-fold with ample canonical divisor.
Then
$$K_X^3\ge \frac{2}{3}(2p_g(X)-7).$$
\endproclaim

When $X$ is not necessarily smooth, we have the following
\proclaim{Theorem 3} Let $X$ be a minimal projective 3-fold
of general type (with only ${\Bbb Q}$-factorial terminal singularities).
Assume that $n+1=p_g(X)\ge 2$ and let $\phi_1:X\dashrightarrow{\Bbb P}^{n}$ be the canonical map.
Then we have the following inequalities according to the dimension of $\phi_1(X):$

(1) $K_X^3\ge 2p_g(X)-6$ if $\dim\phi_1(X)=3.$

(2) $K_X^3\ge p_g(X)-2$ if $\dim\phi_1(X)=2$ and $p_g(X)\ge 6$. If, in addition,
a general fibre of $\phi_1$ is a curve of genus$\ge 3$, then
$K_X^3\ge 2p_g(X)-4.$

(3) When $\phi_1$ is a curve, let $S$ be the minimal model of a general irreducible member of the
movable part of $|K_X|$ and put
$a=K_S^2$,  $b=p_g(S)$. Assume $k=[\frac{1}{2}(p_g-2)]\ge 4$, where $[x]$ stands for the round down of $x$.
Then we have
$$K_X^3\geq
\cases
\text{min}\{\frac{6k^2}{3k^2+8k+4}\cdot (p_g(X)-\frac{4}{3}),
\frac{6k}{3k+4}\cdot (p_g(X)-\frac{5}{3})\},\ &\text{if}\ (a,b)=(1,1)\\
\frac{k^2}{(k+1)^2}\cdot a\cdot (p_g(X)-1), \ &\text{if}\ (a,b)\ne (1,1)
\endcases$$
\endproclaim

The intersection numbers between Weil divisors on singular surfaces are not necessarily
integers, which causes difficulties to get optimal estimates in case (3).

\remark{Remark 4} We make extra assumptions on $p_g(X)$ in Theorem 3(2), 3(3) simply for getting
better inequalities. Our method works also for the case $p_g(X)\ge 2$.
Recall that the geometric genus of a surface of general type with $K_S^2=1$ is bounded by 2 from
above. Furthermore, the surface in case (3) of the theorem has positive geometric genus. Hence Theorem
3 asserts that $K_X^3\ge 2p_g(X)-6$ unless $X$ is canonically fibred by curves of genus two in case (2)
or by surfaces with $a=K_S^2=1$, $b=p_g(S)=2$ in case (3).
\endremark
\medskip

When $X$ is Gorenstein, we have the following theorem, which improves the results known so far:

\proclaim{Theorem 5} Let $X$ be a minimal projective Gorenstein  3-fold of general type with only
locally factorial terminal singularities.

(1) Assume that $X$ is neither canonically fibred by surfaces $S$ with
$c_1(S)^2=1$, $p_g(S)=2$ nor by curves of genus two. Then
$K_X^3\ge 2p_g(X)-6.$

(2) Assume that $X$ is smooth and that $X$ is not canonically fibred by surfaces $S$ with $c_1(S)^2=1$, $p_g(S)=2$. Then
$K_X^3\ge \frac{2}{3}(2p_g(X)-5).$

(3) Assume that the canonical model of $X$ is factorial. If
$K_X^3<\frac{2}{21}(11p_g(X)-16)$, then $X$ is not smooth and is canonically fibred by curves of genus
two.
\endproclaim

These inequalities have a certain interesting application which will be presented in another note.

\head {\rm 1. Preliminaries} \endhead
\subhead 1.1 Conventions \endsubhead
Let $X$ be a normal projective variety of dimension $d$. We denote by
$\Div(X)$
the group of Weil divisors on $X$. An element $D\in\Div(X)\otimes{\Bbb Q}$
is called a ${\Bbb Q}$-{\it divisor}. A ${\Bbb Q}$-divisor $D$ is
said to be ${\Bbb Q}$-{\it Cartier} if $mD$ is a Cartier divisor for some
positive integer $m$. For a ${\Bbb Q}$-Cartier divisor
$D$ and an irreducible curve $C\subset X$, we can define the intersection
number $D\cdot C$ in a natural way. A $\Bbbq$-Cartier divisor $D$ is called
{\it nef} (or {\it numerically effective}) if $D\cdot C\ge 0$ for any
effective curve
$C \subset X$. A nef divisor $D$ is called {\it big} if $D^d>0$. We say that
$X$ is $\Bbbq$-{\it factorial} if every Weil divisor on $X$ is
$\Bbbq$-Cartier.
For a Weil divisor $D$ on $X$, denote by $\Co{X}(D)$ the corresponding
reflexive sheaf.
Denote by $K_X$ a canonical divisor of $X$, which is a Weil divisor. $X$ is
called {\it minimal} if $K_X$ is a nef $\Bbbq$-Cartier divisor. $X$ is said to be of general
type if $\kappa(X)=\dim(X)$. We refer to \cite{R1} for definitions of canonical and terminal singularities.

The symbols $\sim$, $\equiv$ and
$=_{\Bbb Q}$ respectively stands for linear, numerical and
 ${\Bbb Q}$-linear equivalences.

\subhead 1.2 Vanishing theorem \endsubhead
Let $D=\sum a_iD_i$ be a $\Bbbq$-divisor on $X$, where the $D_i$ are distinct
prime divisors and $a_i\in\Bbbq$. We define
$$\align
&\text{the round-down}\ \rounddown{D}:=\sum\rounddown{a_i}D_i,\
\text{where}\
\rounddown{a_i}\ \text{is the integral part of}\ a_i;\\
&\text{the round-up}\  \roundup{D}:=-\rounddown{-D};\\
&\text{the fractional part}\ \{D\}:=D-\rounddown{D}.
\endalign$$
We always use the Kawamata-Viehweg vanishing theorem in the
following form.
\proclaim{Vanishing Theorem} (\cite{Ka1} or \cite{V1}) Let $X$ be a smooth
complete variety, $D\in\Div(X)\otimes\Bbbq$. Assume the following two conditions:

(i) $D$ is nef and big;

(ii) the fractional part of $D$ has supports with only normal crossings.

\noindent Then $H^i(X, \Co{X}(K_X+\roundup{D}))=0$ for all $i>0$.
\endproclaim

Note that, when $S$ is a surface, the above theorem is true without the condition (ii) according to Sakai (\cite{S}) or Miyaoka (\cite{M3, Proposition 2.3}) (also cited in
\cite{E-L, (1.2)}).

\subhead 1.3 Set up for canonical maps \endsubhead
Let $X$ be a projective minimal 3-fold with only ${\Bbb Q}$-factorial terminal singularities.
Suppose $p_g(X)\ge 2$. We study the canonical map $\fei{1}$ which is usually a rational map.
 Take the birational modification $\pi: X'\lrw X$, following Hironaka, such that

(1) $X'$ is smooth;

(2) the movable part of $|K_{X'}|$ is base point free;

(3) $\pi^*(K_X)$ is linearly equivalent to a divisor supported by a divisor of normal crossings.

Denote by $g$ the composition $\fei{1}\circ\pi$. So
$g: X'\lrw W'\subseteq{\Bbb P}^{p_g(X)-1}$
is a morphism. Let
$g: X'\overset{f}\to\lrw W\overset{s}\to\lrw W'$
be the Stein factorization of $g$. We can write
$$K_{X'}=_{\Bbb Q}\pi^*(K_X)+E=_{\Bbb Q}S_1+Z_1,$$
where $S_1$ is the movable part of $|K_{X'}|$, $Z_1$ the fixed part and $E$ is an effective
${\Bbb Q}$-divisor which is a ${\Bbb Q}$-linear combination of distinct exceptional divisors.
We can also write
$$\pi^*(K_X)=_{\Bbb Q} S_1+E',$$
where $E'=Z_1-E$ is actually an effective ${\Bbb Q}$-divisor  and so $\roundup{\pi^*(K_X)}$ means
$\roundup{S_1+E'}$. We note that $1\le \dim (W)\le 3$.

If $\dim\fei{1}(X)=2$, we see that a general fiber of $f$ is a smooth projective curve of genus $g\ge 2$. We say that $X$ is {\it canonically fibred by curves of genus $g$}.

If $\dim\fei{1}(X)=1$, we see that a general fiber $F$ of $f$ is a smooth projective
surface of general type. We say that $X$ is {\it canonically fibred by surfaces}
with invariants $(c_1^2, p_g):=(K_{F_0}^2, p_g(F))$, where $F_0$ is the minimal model of $F$.

\head {\rm 2.  Several simple lemmas} \endhead
The following result is a direct application of an inequality on curves proved by Castelnuovo (\cite{Cas}) and Beauville
(\cite{Be1}).

\proclaim{Lemma 2.1} (\cite{Ch1, Proposition 2.1}) Let $S$ be a smooth projective algebraic surface and $L$ an effective, nef and prime divisor on $S$. Suppose $(K_S-L)\cdot L\ge 0$ and $|L|$ defines a birational rational map onto its image. Then
$$L^2\ge 3h^0(S, \Co{S}(L))-7.$$
\endproclaim

\proclaim{Lemma 2.2} Let $S$ be a smooth projective surface of general type and $L$  a nef divisor on $S$. The following holds.

(i) Suppose that $|L|$ gives a non-birational,  generically finite map onto its image. Then  $L^2\ge 2h^0(S,\Co{S}(L))-4.$

(ii) Suppose that there exists a linear subsystem $\Lambda\subset |L|$ such that $\Lambda$ defines a generically finite map of degree $d$ onto its image. Then
$L^2\ge d[\dim_{\Bbb C}\Lambda-1]$ where $\dim_{\Bbb C}\Lambda$ denotes the projective dimension of $\Lambda$.
\endproclaim
\demo{Proof}
(i) is a special case of (ii).

In order to prove (ii), we take blow-ups  $\pi: S'\lrw S$ such that $\Phi_{\pi^*\Lambda}$ gives a morphism. Let $M$ be the movable part of $\pi^*\Lambda$. Then $h^0(S', M)=\dim_{\Bbb C}\Lambda+1$ and
$$M^2\ge d(h^0(S', M)-2).$$
Since $M\le \pi^*(L)$, we get the inequality
$L^2\ge M^2\ge d(\dim_{\Bbb C}\Lambda-1).$
\qed\enddemo

\proclaim{Lemma 2.3} Let $C$ be a complete smooth algebraic curve. Suppose $D$ is a divisor on $C$ such that $h^0(C, \Co{C}(D))\ge g(C)+1$. Then $\deg(D)\ge 2g(C)$.
\endproclaim
\demo{Proof}
This is a direct result by virtue of R-R and Clifford's theorem.
\qed\enddemo

\proclaim{Lemma 2.4}
Let $S$ be a smooth minimal projective surface of general type. The following holds:

(i) $|mK_S|$ is base point free for all $m\ge 4$;

(ii) $|3K_S|$ is base point free provided $K_S^2\ge 2$;

(iii) $|3K_S|$ is base point free provided $p_g(S)>0$ and $p_g(S)\ne 2$;

(iv) $|2K_S|$ is base point free provided $p_g(S)>0$ or $K_S^2\ge 5$.
\endproclaim
\demo{Proof}
Both (i) and (ii) can be derived from results of Bombieri (\cite{Bo}) and Reider (\cite{Rr}).

If $p_g(S)\ge 3$, then $K_S^2\ge 2$ by Noether inequality. The base point freeness of $|3K_S|$ follows from (ii). If
$K_S^2=1$ and $p_g(S)=1$, $|3K_S|$ is base point free by \cite{Cat1}. If $K_S^2=1$ and $p_g(S)=2$,
 $|3K_S|$ definitely has base points. So (iii) is true.

(iv) follows from  \cite{Ci, Theorem 3.1} and Reider's theorem.
\qed\enddemo

\proclaim{Lemma 2.5} Let $S$ be a smooth projective surface of general type. Let $\sigma: S\lrw S_0$ be the contraction onto the minimal
model. Suppose that there is an effective irreducible curve $C$ on $S$ such that $C\le \sigma^*(2K_{S_0})$
and $h^0(S, C)=2$. If $K_{S_0}^2=p_g(S)=1$, then $C\cdot \sigma^*(K_{S_0})\ge 2$.
\endproclaim
\demo{Proof}
We may assume that $|C|$ is a free pencil. Otherwise, we blow-up $S$ at base points of $|C|$.
Denote $C_1:=\sigma(C)$. Then $h^0(S_0, C_1)\ge 2$. Suppose $C\cdot \sigma^*(K_{S_0})=1$.
Then $C_1\cdot K_{S_0}=1$. Because $p_a(C_1)\ge 2$, we see that $C_1^2>0$.
{}From $K_{S_0}(K_{S_0}-C_1)=0$, we get $(K_{S_0}-C_1)^2\le 0$, i.e. $C_1^2\le 1$. Thus $C_1^2=1$ and
$K_{S_0}\equiv C_1$. This means $K_{S_0}\sim C_1$ by virtue of \cite{Cat1}, which is impossible because
$p_g(S)=1$.  So $C\cdot\sigma^*(K_{S_0})\ge 2$.
\qed
\enddemo

\proclaim{Lemma 2.6} (\cite{Ch4, Lemma 2.7})
Let $X$ be a smooth projective variety of dimension $\ge 2$. Let $D$ be a divisor on $X$ such that
$h^0(X, \Co{X}(D))\ge 2$. Let $S$ be a smooth prime divisor on $X$ and assume that $S$ is not
contained in the fixed part of $|D|$. Denote by $M$ the movable part of $|D|$ and by $N$ the movable
part of $|D|_S|$ on $S$. If the natural restriction map
$$H^0(X, \Co{X}(D))\overset\theta\to\lrw H^0(S, \Co{S}(D|_S))$$
is surjective, then $M|_S\ge N$ and,  in particular,
$$h^0(S,\Co{S}(M|_S))=h^0(S,\Co{S}(N))=h^0(S,\Co{S}(D|_S)).$$
\endproclaim

\head {\rm 3. Proof of Theorem 3} \endhead
We give estimates of $K_X^3$ according to the dimension of the canonical image
$\fei{1}(X)$. Let the notation be as in (1.3) throughout this section. Thus $S_1$ is a general member of the movable part of $|\pi^*(K_X)|$ on a resolution of the indeterminacy of $\phi_1$.
\medskip

{\bf The first case} is $\dim\fei{1}(X)=3$.  Kobayashi (\cite{Kob}) proved

\proclaim{Proposition 3.1}
Let $X$ be a projective minimal algebraic 3-fold of general type with only ${\Bbb Q}$-factorial terminal singularities. Suppose $\dim\fei{1}(X)=3$. Then
$$K_X^3\ge 2p_g(X)-6.$$
\endproclaim
\demo{Proof} We give a very simple proof of this result in order to keep this note self-contained.

In this situation, a general member $S_1\in |S_1|$ is a smooth irreducible projective surface of general type. Because $K_X$ is nef and big, we have
$K_X^3=\pi^*(K_X)^3\ge S_1^3.$
Denote $L:=S_1|_{S_1}$. Then $L$ is a nef and big divisor on $S_1$ and $|L|$ defines a generically
finite map onto its image. It is obvious that
$$h^0(S_1, L)\ge h^0(X', S_1)-1=p_g(X)-1.$$
Note also that $p_g(X)\ge 4$ under the assumption of this proposition.

If $|L|$ gives a birational map, then, by Lemma 2.1,
$$L^2\ge 3h^0(S_1, L)-7
\ge 3p_g(X)-10\ge 2p_g(X)-6.$$

If $|L|$ gives a non-birational rational map, then, by Lemma 2.2,
$$L^2\ge 2h^0(S_1, L)-4\ge 2p_g(X)-6.$$
Therefore
$K_X^3\ge S_1^3=L^2\ge 2p_g(X)-6.$
The proof is complete.
\qed\enddemo

{\bf The second case} is $\dim\fei{1}(X)=2$. The general member $S_1$ is an irreducible smooth surface of general type. The canonical map gives  a fibration
$f: X'\lrw W$, and we  let $C$ denote its general fiber, which is a smooth curve of genus $\ge 2$.

\proclaim{Proposition 3.2}
Let $X$ be a  projective minimal algebraic 3-fold of general type with only ${\Bbb Q}$-factorial terminal singularities. Suppose $\dim\fei{1}(X)=2$ and $p_g(X)\ge 6$. Then either $g(C)\ge 3$ and
$K_X^3\ge \frac{2}{3}g(C)(p_g(X)-2)$ or
$C$ is a curve of genus $2$ and
$K_X^3\ge p_g(X)-2.$
\endproclaim
\demo{Proof} We prove the proposition through several steps.
\medskip

Step 1 (bounding $K_X^3$ in terms of $(L_1, C)$).
Recall that we have $\pi^*(K_X)=_{\Bbb Q}S_1+E'$, where $E'$ is an effective ${\Bbb Q}$-divisor. Put
$L_1:=\pi^*(K_X)|_{S_1}\ \text{and}\ L:=S_1|_{S_1}.$
Then $L_1$ is a nef and big ${\Bbb Q}$-divisor on the surface $S_1$ and $|L|$ is composed of a free pencil of curves on $S_1$. It is obvious that $L_1^2\ge L_1\cdot L$. We can write
$$L=S_1|_{S_1}\sim \sum_{i=1}^a C_i\equiv aC,$$
where
$a\ge h^0(S_1, L)-1\ge p_g(X)-2$
and the $C_i's$ are fibers of $f$ contained in the surface $S_1$. Thus we see that
$$K_X^3=\pi^*(K_X)^3\ge L_1^2\ge L_1\cdot L
\ge (L_1\cdot C)\cdot (p_g(X)-2),$$
and we get a lower bound of $K_X^3$ by giving an estimate of $(L_1\cdot C)$ from
below.
\medskip

Step 2 (the generic finiteness of the tricanonical map $\fei{3}$).
Look at the sublinear system
$$|K_{X'}+\roundup{\pi^*(K_X)}+S_1|\subset |3K_{X'}|.$$
We claim that $\fei{3}$ is generically finite whenever $p_g(X)\ge 4$.
We only have to prove that $\fei{3}|_{S_1}$ is generically finite for a general member $S_1$. By the vanishing theorem, we have
$$\align
|K_{X'}+\roundup{\pi^*(K_X)}+S_1|\bigm|_{S_1}&=
\bigm|K_{S_1}+\roundup{\pi^*(K_X)}|_{S_1}\bigm|\\
&\supset \bigm|K_{S_1}+\roundup{\pi^*(K_X)|_{S_1}}\bigm|.
\endalign$$
We want to prove that
$\Fi{K_{S_1}+\roundup{\pi^*(K_X)|_{S_1}}}$
is generically finite. Because
$ K_{S_1}+\roundup{\pi^*(K_X)|_{S_1}}\ge L,$
we see that $|K_{S_1}+\roundup{\pi^*(K_X)|_{S_1}}|$ separates different fibers of $\Fi{L}$. So we only have to verify that
$\Fi{K_{S_1}+\roundup{\pi^*(K_X)|_{S_1}}}|_C $
is finite for an arbitrary smooth fiber $C$ of $f$ contained in $S_1$.
We have
$$L_1\equiv L+E_{\Bbb Q}\equiv aC+E_{\Bbb Q},$$
where $a\ge p_g(X)-2\ge 2$ and $E_{\Bbb Q}:=E'|_{S_1}$ is an effective ${\Bbb Q}$-divisor on $S_1$.
Thus
$$L_1-C-\frac{1}{a}E_{\Bbb Q}\equiv (1-\frac{1}{a})L_1$$
is a nef and big ${\Bbb Q}$-divisor. Using the vanishing theorem again, we get
$$H^1(S_1, K_{S_1}+\roundup{L_1-\frac{1}{a}E_{\Bbb Q}}-C)=0.$$
This means that
$|K_{S_1}+\roundup{L_1-\frac{1}{a}E_{\Bbb Q}}|\bigm|_C=|K_C+D|,$
where $D:=\roundup{L_1-\frac{1}{a}E_{\Bbb Q}}|_C$ is a divisor on $C$ with positive degree. Because $g(C)\ge 2$, the linear system $|K_C+D|$ gives a finite map, implying the generic finiteness of $\fei{3}$.
\medskip

Step 3 (Estimation of $(L_1\cdot C)$).
Since $|3K_{X'}|$ gives a generically finite map, so does
$\bigm|M_3|_{S_1}\bigm|$ on the surface $S_1$, where $M_3$ is the movable part of $|3K_{X'}|$. Thus $\Fi{M_3|_{S_1}}$ maps general $C$ of genus $\ge 2$ to a curve and hence
$M_3|_{S_1}\cdot C\ge 2.$
Noting that
$3\pi^*(K_X)=_{\Bbb Q} M_3+E_3$
where $E_3$ is an effective ${\Bbb Q}$-divisor,
we see that
$$3\pi^*(K_X)|_{S_1}\cdot C\ge M_3|_{S_1}\cdot C\ge 2,$$
{\it i.e.}, $L_1\cdot C\ge \frac{2}{3}$.
{}From this crude initial estimate, we derive a better one.
To do this, we run a recursive program (the $\alpha$-{\it program}) below.

Pick up a positive integer  $\alpha$. We have
$$|K_{X'}+\roundup{\alpha\pi^*(K_X)}+S_1|\subset |(\alpha+2)K_{X'}|.$$
The vanishing theorem gives
$$\align
|K_{X'}+\roundup{\alpha\pi^*(K_X)}+S_1|\bigm|_{S_1}&=
\bigm|K_{S_1}+\roundup{\alpha\pi^*(K_X)}|_{S_1}\bigm|\\
&\supset |K_{S_1}+\roundup{\alpha L_1}|.
\endalign$$
We see that
$\alpha L_1-C-\frac{1}{a}E_{\Bbb Q}\equiv (\alpha-\frac{1}{a})L_1$
is a nef and big ${\Bbb Q}$-divisor. Using the vanishing theorem on $S_1$ again, we get
$$|K_{S_1}+\roundup{\alpha L_1-\frac{1}{a}E_{\Bbb Q}}|\bigm|_C=
|K_C+D_{\alpha}|, \tag{3.1}$$
where $D_{\alpha}:=\roundup{\alpha L_1-\frac{1}{a}E_{\Bbb Q}}|_C$ with
$\deg(D_{\alpha})\ge \roundup{(\alpha-\frac{1}{a})L_1\cdot C}.$
We have to use several symbols in order to obtain our result. Let $M_{\alpha+2}$
be the movable part of $|(\alpha+2)K_{X'}|$. Let $M_{\alpha+2}'$ be the movable part of
$$|K_{X'}+\roundup{\alpha\pi^*(K_X)}+S_1|.$$
Clearly we have $M_{\alpha+2}'\le M_{\alpha+2}.$
Let $N_{\alpha}$ be the movable part of $|K_{S_1}+\roundup{\alpha L_1}|$. Then it is easy to see $M_{\alpha+2}'|_{S_1}\ge N_{\alpha}$ by Lemma 2.6. So
$$(\alpha+2)L_1\ge_{\Bbb Q}M_{\alpha+2}|_{S_1}\ge M_{\alpha+2}'|_{S_1}
\ge N_{\alpha}.$$
Let $N_{\alpha}'$ be the movable part of
$|K_{S_1}+\roundup{\alpha L_1-\frac{1}{a}E_{\Bbb Q}}|.$
Then obviously $N_{\alpha}\ge N_{\alpha}'$.
{}From (3.1) and Lemma 2.6, we have $h^0(C, N_{\alpha}'|_C)=h^0(C, K_C+D_{\alpha}).$
Thus we see that
$$h^0(C, N_{\alpha}|_C)\ge h^0(C, N_{\alpha}'|_C)=h^0(C, K_C+D_{\alpha}).$$

Now take $\alpha=2$ and run the $\alpha$-program. We get
$4L_1\cdot C\ge N_2\cdot C.$
Because $a>3$ under the assumption, we see that
$\deg(D_2)\ge \roundup{(2-\frac{1}{a})\frac{2}{3}}= 2.$
Thus $h^0(C, N_2|_C)\ge g(C)+1$. By Lemma 2.3, we have
$N_2\cdot C\ge 2g(C).$
If $g(C)=2$, we get $L_1\cdot C\ge 1$ and thus the inequality
$K_X^3\ge p_g(X)-2.$
If $g(C)\ge 3$, we get
$L_1\cdot C\ge \frac{3}{2}.$
This is a better bound than the initial one. However this is not enough to derive our statement. We have to optimize our estimation.
\medskip

Step 4 (Optimization).
As has been seen in the previous step, we have $L_1\cdot C\ge \frac{3}{2}$ when $g\ge 3$.
We take $\alpha=1$ now and run the $\alpha$-program. Since $p_g(X)\ge 6$, we have $a\ge 4$. Thus
$$\deg(D_1)\ge\roundup{(1-\frac{1}{a})\frac{3}{2}}= 2.$$
So $h^0(C, N_1|_C)\ge g(C)+1$. Therefore we get, by Lemma 2.3, that
$$3L_1\cdot C\ge N_1\cdot C\ge 2g(C)\ge 6 \ \text{whenever}\ g(C)\ge 3.$$
This means $L_1\cdot C\ge 2$, which is what we want. So we have the inequality
$$K_X^3\ge \frac{2}{3}\cdot g(C)\cdot (p_g(X)-2)\tag{3.2}$$
whenever $g(C)\ge 3$. The proof is complete.
\qed
\enddemo

{\bf The last case} is $\dim\fei{1}(X)=1$. The canonical map gives a fibration $f:X'\lrw W$ where $W$ is a smooth projective curve. Denote $b:=g(W)$. We see that a general fiber $F$ of $f$ is a smooth projective surface of general type. Let $\sigma: F\lrw F_0$ be the contraction onto the minimal model. Note that we always have $p_g(F)>0$ in this situation. We also have
$S_1\sim \sum_{i=1}^{b_1} F_i\equiv {b_1}F,$
where the $F_i's$ are fibers of $f$ and $b_1\ge p_g(X)-1$.

\proclaim{Proposition 3.3}
Let $X$ be a projective minimal algebraic 3-fold of general type with only ${\Bbb Q}$-factorial terminal singularities. Suppose $\dim\fei{1}(X)=1$. Let $k\ge 4$ be an integer and assume that $p_g(X)\ge 2k+2$. Then
$K_X^3\ge\frac{k^2}{(k+1)^2}\cdot K_{F_0}^2\cdot (p_g(X)-1).$
\endproclaim
\demo{Proof} The proof proceeds through two steps.
\medskip

Step 1 (bounding $K_X^3$ in terms of $L^2$).
On the surface $F$, we denote $L:=\pi^*(K_X)|_F$. Then $L$ is an effective nef and big ${\Bbb Q}$-divisor. Because
$\pi^*(K_X)\equiv b_1F+E'$
with $E'$ effective, we get
$$K_X^3=\pi^*(K_X)^3\ge (\pi^*(K_X)^2\cdot F)\cdot (p_g(X)-1)=L^2\cdot (p_g(X)-1).$$
So the main point is to estimate $L^2$ from below in order to prove the proposition.
\medskip

Step 2 (bounding $L^2$ from below by studying the $(k+1)$-canonical map $\fei{k+1}$).
Let $M_{k+1}$ be the movable part of $|(k+1)K_{X'}|$. Then we may write
$$(k+1)\pi^*(K_X)=_{\Bbb Q}M_{k+1}+E_{k+1}$$
where $E_{k+1}$ is an effective ${\Bbb Q}$-divisor. Therefore we see that
$(k+1)L\ge_{\text{num}} M_{k+1}|_F$. Let $N_k$ be the movable part of $|kK_F|$. According to Lemma 2.4, $|kK_{F_0}|$ is base point free.  Thus $N_k=\sigma^*(kK_{F_0})$. We claim that $M_{k+1}|_F\ge N_k$. Then $(k+1)L\ge N_k$ and we get
$$L^2\ge\frac{1}{(k+1)^2}N_k^2=\frac{k^2}{(k+1)^2}K_{F_0}^2.$$
So we have the inequality
$$K_X^3\ge\frac{k^2}{(k+1)^2}\cdot K_{F_0}^2\cdot (p_g(X)-1). \tag{3.3}$$

Now we prove the claim. In fact, $\fei{1}$ is a morphism if $b>0$. In this case, we
do not need any modification and $f:X'=X\lrw W$ is a fibration. A general fiber $F$ is a smooth projective surface of general type, because the singularities on $X$ are isolated. Furthermore $F$ is minimal because $K_X$ is nef. By Kawamata's vanishing theorem for ${\Bbb Q}$-Cartier Weil divisor (\cite{KMM}), we have
$H^1(X, kK_X)=0$. This means
$|kK_X+F|\bigm|_F=|kK_F|.$
Noting that $F\le K_X$ and using Lemma 2.6, we see that the claim is true in this case.

We then consider the case with $b=0$. We use the approach in \cite{Kol, Corollary 4.8} to prove it. The canonical map gives a fibration $f:X'\lrw {\Bbb P}^1$. Because
$p_g(X)\ge 2k+2$,
we see that ${\Cal O}(2k+1)\hookrightarrow f_*\omega_{X'}$. Thus we have
$${\Cal E}:={\Cal O}(1)\otimes f_*\omega_{X'/{\Bbb P}^1}^k={\Cal O}(2k+1)\otimes f_*\omega_{X'}^k\hookrightarrow f_*\omega_{X'}^{k+1}.$$
Note that
$H^0({\Bbb P}^1, f_*\omega_{X'}^{k+1})\cong H^0(X', \omega_{X'}^{k+1}).$
It is well known that ${\Cal E}$ is generated by global sections and that
$f_*\omega_{X'/{\Bbb P}^1}^k$ is a sum of line bundles with non-negative degree (cf. \cite{F}, \cite{V2, V3}). Thus the global sections of ${\Cal E}$ separates different fibers of $f$. On the other hand, the local sections of $f_*\omega_{X'}^k$ give the k-canonical map of $F$ and these local sections can be extended to global sections of ${\Cal E}$. This essentially means $M_{k+1}|_F\ge N_k$.
\qed
\enddemo

\proclaim{Proposition 3.4}
Let $X$ be a projective minimal algebraic 3-fold of general type with only ${\Bbb Q}$-factorial terminal singularities.
Suppose that $\dim\fei{1}(X)=1$. Let $k\ge 3$ be an integer and assume $p_g(X)\ge 2k+2$.  If $(K_{F_0}^2, p_g(F))=(1,1)$, then
$$K_X^3\ge\text{min}\{\frac{6k^2}{3k^2+8k+4}\cdot (p_g(X)-\frac{4}{3}),
\frac{6k}{3k+4}\cdot (p_g(X)-\frac{5}{3})\}.$$
\endproclaim
\demo{Proof}
{}From Step 2 in the proof of Proposition 3.3, we have shown that
$$(k+1)\pi^*(K_X)|_F\ge M_{k+1}|_F\ge k\sigma^*(K_{F_0}).$$
(Although we suppose $k\ge 4$ in Proposition 3.3, the case with $k=3$ can be parallelly treated since $|3K_{F_0}|$ is base point free for a surface with
$(K_{F_0}^2, p_g(F))=(1,1)$.)

The canonical map derives a fibration $f: X'\lrw W$. Because $q(F)=0$, we have
$$q(X)=h^1(\Co{X'})=b+h^1(W, R^1f_*\omega_{X'})=b,$$
$$\align
h^2(\Co{X})&=h^1(W, f_*\omega_{X'})+h^0(W, R^1f_*\omega_{X'})\\
&=h^1(W, f_*\omega_{X'})\le 1.
\endalign$$

It is obvious that $h^2(\Co{X})=0$ when $b=0$, since $f_*\omega_{X'}$ is a line bundle of positive degree. Anyway, we have $q(X)-h^2(\Co{X})\ge 0$.
Thus we get
$$\chi(\omega_{X})=p_g(X)+q(X)-h^2(\Co{X})-1\ge p_g(X)-1.$$
By the plurigenus formula of Reid (\cite{R1}), we have
$$P_2(X)\ge \frac{1}{2}K_X^3-3\chi(\Co{X})\ge
\frac{1}{2}K_X^3+3[p_g(X)-1]. \tag{3.4}$$
Let $M_2$ be the movable part of $|2K_{X'}|$. We consider the natural restriction map $\gamma$:
$$H^0(X', M_2)\overset{\gamma}\to\lrw V_2\subset H^0(F, M_2|_F)\subset
H^0(F, 2K_F),$$
where $V_2$ is the image of $\gamma$ as a ${\Bbb C}$-subspace of $H^0(F, M_2|_F)$. Because $h^0(2K_F)=3$, we see that
$1\le\dim_{\Bbb C}V_2\le 3$. Denote by $\Lambda_2$ the linear system corresponding to$V_2$. We have $\dim\Lambda_2=\dim_{\Bbb C}V_2-1$.
\medskip

Case 1. $\dim_{\Bbb C}V_2=3$.

Since $\Lambda_2$ is a sub-system of $|2K_F|$, we see that the restriction of $\phi_{2,X'}$ to $F$ is exactly the bicanonical map of $F$. Because $\phi_{2, F}$ is a generically finite morphism of degree $4$, $\phi_{2,X'}$ is also a generically finite map of degree $4$.
Let $S_2\in|M_2|$ be a general member. We can further modify $\pi$ such that
$|M_2|$ is base point free. Then $S_2$ is a smooth projective irreducible surface of general type. On the surface $S_2$, denote $L_2:=S_2|_{S_2}$.
$L_2$ is a nef and big divisor. We have
$$2\pi^*(K_X)|_{S_2}\ge S_2|_{S_2}=L_2.$$
We consider the natural map
$$H^0(X', S_2)\overset{\gamma'}\to\lrw \overline{V_2}\subset H^0(S_2, L_2),$$
where $\overline{V_2}$ is the image of $\gamma'$.
Denote by $\overline{\Lambda_2}$ the linear system corresponding to $\overline{V_2}$.
Because $\fei{2}$ is generically finite map of degree $4$, we see that $|L_2|$ has a sub-system
$\overline{\Lambda_2}$ which gives a generically finite map of degree $4$. By Lemma 2.2(ii), we get
$L_2^2\ge 4(\dim_{\Bbb C}\overline{\Lambda_2}-1)\ge 4(P_2(X)-3).$
Therefore we have
$$K_X^3\ge\frac{1}{8}L_2^2\ge\frac{1}{2}(P_2(X)-3)
\ge \frac{1}{2}(\frac{1}{2}K_X^3+3p_g(X)-6).$$
Therefore
$$K_X^3\ge 2p_g(X)-4. \tag{3.5}$$
\medskip

Case 2. $\dim_{\Bbb C}V_2=2$.

In this case, $\dim\fei{2}(F)=1$ and $\dim\fei{2}(X)=2$. We may further modify
$\pi$ such that $|M_2|$ is base point free. Taking the Stein factorization of
$\fei{2}$, we get a derived fibration $f_2:X'\lrw W_2$ where $W_2$ is
a surface. Let $C$ be a general fiber of $f_2$. we see that $F$ is naturally fibred by curves with the same numerical type as $C$. On the surface $F$, we have a free pencil $\Lambda_2\subset |2K_F|$. Let $|C_0|$ be the movable part of $\Lambda_2$. Then $h^0(F, C_0)=2$. Because $q(F)=0$, we see that $|C_0|$ is a pencil over the rational curve. So
a general member of $|C_0|$ is an irreducible curve. According to Lemma 2.5, we have $(C_0\cdot\sigma^*(K_{F_0}))_F\ge 2$ whence
$$(\pi^*(K_X)\cdot C)_{X'}=(\pi^*(K_X)|_F\cdot C_0)_F\ge\frac{k}{k+1}(\sigma^*(K_{F_0})\cdot C_0)_F\ge \frac{2k}{k+1}.$$
Now we study on the surface $S_2$. We may write
$$S_2|_{S_2}\sim\sum_{i=1}^{a_2}C_i\equiv a_2C,$$
where the $C_i's$ are fibers of $f_2$ and $a_2\ge P_2(X)-2$. Noting that
$$(\pi^*(K_X)|_{S_2}\cdot C)_{S_2}=(\pi^*(K_X)\cdot C)_{X'}\ge\frac{2k}{k+1}$$
and $2\pi^*(K_X)|_{S_2}\ge S_2|_{S_2}$, we get
$$\align
4K_X^3&\ge 2\pi^*(K_X)^2\cdot S_2=2(\pi^*(K_X)|_{S_2})_{S_2}^2\\
  &\ge a_2(\pi^*(K_X)|_{S_2}\cdot C)_{S_2}
\ge\frac{2k}{k+1}(P_2(X)-2)\\
&\ge \frac{2k}{k+1}(\frac{1}{2}K_X^3+3p_g(X)-5).
\endalign$$
Equivalently
$$K_X^3\ge\frac{6k}{3k+4}p_g(X)-\frac{10k}{3k+4}. \tag{3.6}$$
\medskip

Case 3. $\dim_{\Bbb C}V_2=1$.

In this case,  $\dim\fei{2}(X)=1$. Because $p_g(X)>0$, we see that both $\fei{2}$ and
$\fei{1}$ give the same fibration $f:X'\lrw W$ after taking the Stein factorization of them. So we may write
$$
2\pi^*(K_X)\sim \sum_{i=1}^{a_2'}F_i+E_2'
\equiv a_2'F+E_2',$$
where the $F_i's$ are fibers of $f$, $E_2'$ is an effective ${\Bbb Q}$-divisor, $a_2'\ge P_2(X)-1$ and $F$ is a surface with $(K_{F_0}^2, p_g(F))=(1,1).$
So we get
$$\align
2K_X^3&\ge a_2'(\pi^*(K_X)|_F)_F^2\ge\frac{k^2}{(k+1)^2}(P_2(X)-1)\\
&\ge \frac{k^2}{(k+1)^2}(\frac{1}{2}K_X^3+3p_g(X)-4).
\endalign$$
Equivalently
$$K_X^3\ge\frac{6k^2}{3k^2+8k+4}p_g(X)-\frac{8k^2}{3k^2+8k+4}. \tag{3.7}$$

Comparing (3.5), (3.6) and (3.7), we get the inequality.
\qed
\enddemo

Propositions 3.1, 3.2, 3.3 and 3.4 imply Theorem 3.

\head {\rm 4. Inequalities for minimal Gorenstein 3-folds} \endhead
This section is devoted to study lower bounds for $K_X^3$ of Gorenstein 3-folds. Let $X$ be a projective minimal Gorenstein
3-fold of general type with only locally factorial terminal singularities. It is well known that $K_X^3$ is a
positive even integer and $\chi(\Co{X})<0$. We also have the Miyaoka-Yau inequality (\cite{M2}):
$K_X^3\le -72\chi(\Co{X}).$ Besides, after taking a special birational modification to $X$ according to Reid (\cite{R2})
while using a result of Miyaoka (\cite{M2}), we get the plurigenus formula as follows.
$$P_m(X)=(2m-1)(\frac{m(m-1)}{12}K_X^3-\chi(\Co{X})). \tag{4.1}$$

The following theorem improves \cite{Kob, Main Theorem}, where
we use the same notations as in previous sections.

\proclaim{Theorem 4.1}
Let $X$ be a projective minimal Gorenstein 3-fold of general type with only locally factorial terminal singularities.
Then we have

(i) If $\dim\fei{1}(X)=3$, then $K_X^3\ge 2p_g(X)-6$.

(ii) If $\dim\fei{1}(X)=2$, {\it i.e.}, $X$ is canonically fibered by curves of genus $g$,  then
$$K_X^3\ge \roundup{\frac{2}{3}(g-1)}(p_g(X)-2).$$

(iii) If $\dim\fei{1}(X)=1$, then either
$K_X^3\ge 2p_g(X)-4$
or $(K_{F_0}^2, p_g(F))=(1, 2).$
\endproclaim
\demo{Proof}
By Proposition 3.1, it is sufficient to study the cases $\dim\fei{1}(X)<3$.
\medskip

Case 1. $\dim\fei{1}(X)=2$.

The canonical map gives a fibration $f:X'\lrw W$, where a general fiber $C$ is a smooth curve of
genus $g$.
If $g=2$, our inequality is $K_X^3\ge p_g(X)-2$, which is trivially true. Now
we assume $g\ge 3$. Denote $L:=\pi^*(K_X)|_{S_1}$, which is a nef and big Cartier divisor. Let $S_1\in|M_1|$ be a general member. Then $S_1$ is a smooth projective surface of general type. Noting that $|S_1|_{S_1}|$ is composed of a free pencil of curves with the same numerical type as $C$, we have
$$\pi^*(K_X)|_{S_1}\equiv aC+E_2,$$
where $E_2$ is effective and $a\ge p_g(X)-2$, and we immediately see
$$K_X^3\ge (L\cdot C)(p_g(X)-2).$$
Thus it is sufficient to bound $(L\cdot C)$ from below.

We run once more a recursive program (the $\beta$-{\it program}) which is essentially similar to the $\alpha$-{\it program}. There is, however, a minor difference between them.
Pick up a positive integer $\beta$. Obviously, we have
$$|K_{X'}+\beta\pi^*(K_X)+S_1|\subset |(\beta+2)K_{X'}|.$$
The vanishing theorem gives
$$|K_{X'}+\beta\pi^*(K_X)+S_1|\bigm|_{S_1}=|K_{S_1}+\beta L|.$$
We have $L\ge C$. If $\beta>1$, then we have
$$|K_{S_1}+(\beta-1)L+C|\bigm|_C=|K_C+D_{\beta}|,$$
where $D_{\beta}:=(\beta-1)L|_C$.
Let $M_{\beta+2}$ be the movable part of $|(\beta+2)K_{X'}|$ and $M_{\beta+2}'$ be the movable part of $|K_{X'}+\beta\pi^*(K_X)+S_1|$. Then
$M_{\beta+2}\ge M_{\beta+2}'$. Let $N_{\beta}$ be the movable part of $|K_{S_1}+(\beta-1)L+C|$. Then, by Lemma 2.6, we have
$$(\beta+2)L\ge M_{\beta+2}|_{S_1}\ge  M_{\beta+2}'|_{S_1}\ge N_{\beta}.$$
Also by Lemma 2.6, we have $h^0(C, N_{\beta}|_C)=h^0(K_C+D_{\beta})$.
If $\deg(D_{\beta})=(\beta-1)(L\cdot C)\ge 2$,
then
$$h^0(C, N_{\beta}|_C)=g-1+(\beta-1)(L\cdot C).$$
Using R-R again and Clifford's theorem, we see that $h^1(C, N_{\beta}|_C)=0$ and
$$(\beta+2)(L\cdot C)\ge N_{\beta}\cdot C=2g-2+(\beta-1)(L\cdot C).$$
We get the inequality
$$L\cdot C\ge \frac{2g-2+(\beta-1)(L\cdot C)}{\beta+2}. \tag{4.2}$$
Now take $\beta=3$. Then $\deg(D_3)\ge 2$. According to (4.2), we see
$L\cdot C>1$, i.e. $L\cdot C\ge 2$. {}From now on, we can constantly take $\beta=2$. we see that $\deg(D_2)\ge 2$. So (4.2) becomes
$L\cdot C\ge\frac{2g-2}{3}.$
This means $L\cdot C\ge \roundup{\frac{2}{3}(g-1)}.$
\medskip

Case 2. $\dim\fei{1}(X)=1$.

In this case, the canonical map derives a fibration $f:X'\lrw W$ onto a smooth curve $W$ where a general fiber $F$ of $f$ is a smooth irreducible surface of general type. We have $\pi^*(K_X)=S_1+E'$ and $S_1\equiv b_1F$, where $b_1\ge p_g(X)-1$. Denote $\overline{S}=\pi(S_1)$ and $\overline{F}=\pi(F)$. Then
$\overline{S}\equiv b_1\overline{F}$. Because $\overline{F}^2$ is pseudo-effective, $K_X\cdot\overline{F}^2\ge 0$. Note that $K_X\cdot\overline{F}^2$ is an even integer.

If $K_X\cdot\overline{F}^2>0$, then we have $K_X^2\cdot \overline{F}\ge 2(p_g(X)-1)$ and thus
$K_X^3\ge 2(p_g(X)-1)^2.$

If $K_X\cdot\overline{F}^2=0$, then
$\Co{F}(\pi^*(K_X)|_F)\cong\Co{F}(\sigma^*(K_{F_0}))$ by a trivial generalization of \cite{Ch3, Lemma 2.3}.
 Thus we always have
$$\align
K_X^3&=\pi^*(K_X)^3\ge (\pi^*(K_X)^2\cdot F)(p_g(X)-1)\\
&=\sigma^*(K_{F_0})^2(p_g(X)-1)\ge 2(p_g(X)-1)
\endalign$$
whenever $K_{F_0}^2\ge 2$.

When $K_{F_0}^2=1$, the only possibility is $1\le p_g(F)\le 2$. We can prove that
$K_X^3\ge 2p_g(X)-4$
if $(K_{F_0}^2, p_g(F))=(1,1)$. In fact, this is the special case of
Proposition 3.4 and the estimation here is more exact since $X$ is Gorenstein.
The main point is that we have $\pi^*(K_X)|_F\sim \sigma^*(K_{F_0}).$
We see from the proof of Proposition 3.4 that
(3.5) is still as $K_X^3\ge 2p_g(X)-4$,
that (3.6) corresponds to
$K_X^3\ge 2p_g(X)-3\frac{1}{3}$
and that (3.7) will be replaced by
$K_X^3\ge 2p_g(X)-2\frac{2}{3}.$
\qed\enddemo

{}From Theorem 4.1, one sees that bad cases possibly occur when $X$ is
 canonically fibered by curves of genus 2 or by surfaces with invariants $(c_1^2, p_g)=(1,2)$.
For technical reasons, we are only able to treat a nonsingular 3-fold.
One needs a new method to cover singular 3-folds.

Now suppose that $X$ is a smooth projective 3-fold.
Let $\overline{M}$ be a divisor on $X$ such that $h^0(X,\overline{M})\ge 2$ and that $|\overline{M}|$ has base points but no fixed part.
By Hironaka's theorem (\cite{Hi}), we may take successive blow-ups
$$\pi: X'=X_{n}\overset{\pi_n}\to\rw X_{n-1}\rw\cdots\rw
   X_{i}\overset{\pi_i}\to\rw X_{i-1}\rw\cdots\rw
  X_1\overset{\pi_1}\to\rw X_0=X$$
such that

(i) $\pi_i$ is a single blow-up along smooth center $W_i$ on $X_{i-1}$ for all $i$;

(ii) $W_i$ is contained in the base locus of the movable part of
$$|(\pi_1\circ\pi_2\circ\cdots\circ\pi_{i-1})^*(\overline{M})|$$
and thus $W_i$ is a reduced closed point or a smooth projective curve on $X_{i-1}$;

(iii) the movable part of $|\pi^*(\overline{M})|$ has no base points.

It is clear that the resulting 3-fold $X'$ is still smooth.
Let $E_i$ be the exceptional divisor on $X'$ corresponding to $W_i$. Then we may write
$$K_{X'}=\pi^*(K_X)+\sum_{i=1}^n a_iE_i,\ \
\pi^*(\overline{M})=M+\sum_{i=1}^ne_iE_i,$$
where $a_i,\ e_i\in{\Bbb Z}$, $a_i\ge 0$ and $M$ is the movable part of $|\pi^*(\overline{M})|$. {}From the definition of $\pi$, we see $e_i>0$ for all $i$.

\proclaim{Lemma 4.2} $a_i\le 2e_i$ for all $i$.
\endproclaim
\demo{Proof}
We  prove the simple lemma by induction. Denote by $M_i$ the strict transform of $\overline{M}$ in $X_i$ for all $i$. Let $E_i^{(i)}$ be the exceptional divisor
on $X_i$ corresponding to $W_i$. Let $E_i^{(j)}$ be the strict transform of
$E_i^{(i)}$ in $X_j$ for $j>i$.

For $i=1$, we have
$$K_{X_1}=\pi_1^*(K_X)+a_1^{(1)}E_1^{(1)}\ \ \text{and}\ \
\pi_1^*(\overline{M})=M_1+e_1^{(1)}E_1^{(1)}.$$
{}From the definition of $\pi_1$, we know that $e_1^{(1)}\ge 1$. Note that
$a_1^{(1)}$ is computable. In fact,
$a_1^{(1)}=2$ if $W_1$ is a reduced smooth point of $X$;
$a_1^{(1)}=1$ if $W_1$ is a smooth curve on $X$.
Clearly, we have $a_1^{(1)}\le 2e_1^{(1)}$.

For $i=n-1$, we have
$$K_{X_{n-1}}=(\pi_1\circ\cdots\circ\pi_{n-1})^*(K_X)+
\sum_{i=1}^{n-1}a_i^{(n-1)}E_i^{(n-1)}$$
$$(\pi_1\circ\cdots\circ\pi_{n-1})^*(\overline{M})=M_{n-1}+
\sum_{i=1}^{n-1}e_i^{(n-1)}E_i^{(n-1)}.$$
Suppose we have already had $a_i^{(n-1)}\le 2e_i^{(n-1)}$. Then we get
$$\align
K_{X_n}&=\pi_n^*(K_{X_{n-1}})+a_n^{(n)}E_n^{(n)}\\
&=\pi^*(K_X)+\pi_n^*\sum_{i=1}^{n-1}a_i^{(n-1)}E_i^{(n-1)}+a_n^{(n)}E_n^{(n)}.
\endalign$$
$$\align
\pi^*(\overline{M})&=\pi_n^*(M_{n-1})+
\pi_n^*\sum_{i=1}^{n-1}e_i^{(n-1)}E_i^{(n-1)}\\
&=M+\pi_n^*\sum_{i=1}^{n-1}e_i^{(n-1)}E_i^{(n-1)}+e_n^{(n)}E_n^{(n)}.
\endalign$$
Because $\pi_n$ is also a single blow-up, we see similarly that
$a_n^{(n)}\le 2e_n^{(n)}$. Note that  $E_n^{(n)}=E_n$ and
$$\sum_{i=1}^na_iE_i=
\pi_n^*\sum_{i=1}^{n-1}a_i^{(n-1)}E_i^{(n-1)}+a_n^{(n)}E_n;$$
$$\sum_{i=1}^ne_iE_i=
\pi_n^*\sum_{i=1}^{n-1}e_i^{(n-1)}E_i^{(n-1)}+e_n^{(n)}E_n.$$
We see that $a_i\le 2e_i$. The proof is complete.
\qed
\enddemo

\proclaim{Theorem 4.3}
Let $X$ be a projective minimal smooth 3-fold of general type.
Suppose $\dim\fei{1}(X)=2$ and $X$ is canonically fibred by curves of genus $2$. Then
$$K_X^3\ge \frac{1}{3}(4p_g(X)-10).$$
The inequality is sharp.
\endproclaim
\demo{Proof}
We keep the same notations as in 1.3 and in Case 1 of the proof of Theorem 4.1.
Set $K_X\sim \overline{M}+\overline{Z}$, where $\overline{M}$ is the movable part of $|K_X|$ and $\overline{Z}$ is the fixed part.
We may take the same successive blow-ups
$$\pi: X'=X_{n}\overset{\pi_n}\to\rw X_{n-1}\rw\cdots\rw
   X_{i}\overset{\pi_i}\to\rw X_{i-1}\rw\cdots\rw
  X_1\overset{\pi_1}\to\rw X_0=X$$
as in the set up for Lemma 4.2.

Let $g=\fei{1}\circ\pi$. Taking the Stein-factorization of $g$, we get the induced fibration $f:X'\lrw W$. A general fiber of $f$ is a smooth curve of genus $2$ by assumption of the theorem. Let $S_1$ be the movable part of
$|\pi^*(\overline{M})|$. Then we have
$$K_{X'}=\pi^*(K_X)+E=\pi^*(K_X)+\sum_{i=0}^pa_iE_i$$
and $\pi^*(\overline{M})\sim S_1+\sum_{i=0}^pe_iE_i.$
We know that $a_i\ge 0$, $e_i>0$ and both $a_i$ and $e_i$ are integers for all $i$. We also have
$$\align
\pi^*(K_X)&=\pi^*(\overline{M})+\pi^*(\overline{Z})=S_1+\sum_{i=0}^pe_iE_i+
\pi^*(\overline{Z})\\
&\sim S_1+\sum_{i=0}^pe_i'E_i+\sum_{j=1}^q d_jL_j
=S_1+E',\endalign$$
where $e_i'\ge e_i$, $d_j>0$, $E_i\ne L_j$ and $L_{j_1}\ne L_{j_2}$ provided
$j_1\ne j_2$.
On the surface $S_1$, set $L:=\pi^*(K_X)|_{S_1}$. We also have
$S_1|_{S_1}\equiv aC$ where $a\ge p_g(X)-2$ and $C$ is a general fiber of the restricted fibration
$f|_{S_1}: S_1\lrw f(S_1).$
Note that the above $C$ lies in the same numerical class as that of a general fiber of $f$. If $L\cdot C\ge 2$, we have already seen in the proof of Theorem 4.1 that
$K_X^3\ge 2p_g(X)-4.$
{}From now on, we suppose $L\cdot C=1$. Note that, in this situation, $|\overline{M}|$ definitely has base points. (Otherwise, $\pi=$identity and
$$L\cdot C=K_X|_{S_1}\cdot C=(K_X+S_1)|_{S_1}\cdot C=K_{S_1}\cdot C=2$$
which contradicts to the assumption $L\cdot C=1$.)

Denote $E'|_{S_1}:=E_V'+E_H'$, where $E_V'$ is the vertical part, {\it i.e.,}
$\dim f|_{S_1}(E_V')=0$, and $E_H'$ is the horizontal part, {\it i.e.,} $E_H'\cdot C>0$.
Because $E'|_{S_1}\cdot C=L\cdot C=1$, we see that $E_H'\cdot C=1$. This means that $E_H'$ is an irreducible curve and is a section of the restricted fibration $f|_{S_1}$. Denote $E|_{S_1}:=E_V+E_H$, where $E_V$ is the vertical part and $E_H$ is the horizontal part. {}From $K_{S_1}\cdot C=2$, one sees that
$E_H\cdot C=E|_{S_1}\cdot C=1$. This also means that $E_H$ is an irreducible curve and $E_H$ comes from some exceptional divisor $E_i$ with $a_i=1$. We may suppose that $E_H$ comes from $E_0$. Then $a_0=1$. Because $e_0'>0$ and $\pi^*(K_X)\cdot C=1$, we see that $e_0'=1$ and thus $E_H'$ also comes from $E_0$. Since $E_0|_{S_1}$ has only one horizontal part, $E_H$ and $E_H'$ coincide with a curve $G$. Now It is quite clear that
$$E_V=\sum_{i=1}^pa_i(E_i|_{S_1})+(E_0|_{S_1}-G),$$
$$E_V'=\sum_{i=1}^pe_i'(E_i|_{S_1})+\sum_{j=1}^qd_j(L_j|_{S_1})
+(E_0|_{S_1}-G).$$
We have the following
\medskip

\noindent{\bf Claim.}\ \ $E_V\le 2E_V'.$
\medskip

This is apparently a direct consequence of Lemma 4.2. In fact, we have $a_i\le 2e_i\le 2e_i'$ by Lemma 4.2 for all $i>0$. Thus
$$\sum_{i=1}^pa_i(E_i|_{S_1})\le 2\sum_{i=1}^pe_i'(E_i|_{S_1})
\le 2(\sum_{i=1}^pe_i'(E_i|_{S_1})+\sum_{j=1}^qd_j(L_j|_{S_1})).$$
On the other hand, It is obvious that
$E_0|_{S_1}-G\le 2(E_0|_{S_1}-G).$
Therefore we get
$$\align
E_V&=(E_0|_{S_1}-G)+\sum_{i=1}^pa_i(E_i|_{S_1})\\
&\le 2(E_0|_{S_1}-G)
+2(\sum_{i=1}^pe_i'(E_i|_{S_1})+\sum_{j=1}^qd_j(L_j|_{S_1}))
=2E_V'
\endalign$$
and the claim is proved.

Since that $2E_V'-E_V$ is effective and vertical, we see that
$E_V\cdot G\le 2E_V'\cdot G$. On the surface $S_1$, we have
$$(K_{S_1}+2C+G)G=2p_a(G)-2+2G\cdot C=2p_a(G)\ge 0.$$
On the other hand, we have
$$\align
&(K_{S_1}+2C+G)G\\
=&((\pi^*(K_X)|_{S_1}+E_V+G+S_1|_{S_1})+2C+G)G\\
\le & (\pi^*(K_X)|_{S_1}+S_1|_{S_1}+G)\cdot G+2E_V'\cdot G+2+G^2\\
=& 2\pi^*(K_X)|_{S_1}\cdot G+E_V'\cdot G+G^2+2.
\endalign$$
So we have
$$2\pi^*(K_X)|_{S_1}\cdot G+E_V'\cdot G+G^2+2\ge 0. \tag{4.3}$$
We also have
$$\pi^*(K_X)|_{S_1}\cdot G=S_1|_{S_1}\cdot G+E_V'\cdot G+G^2. \tag{4.4}$$
Combining (4.3) and (4.4), we get
$$3\pi^*(K_X)|_{S_1}\cdot G\ge S_1|_{S_1}\cdot G-2\ge p_g(X)-4.$$
$$\pi^*(K_X)\cdot S_1\cdot E'\ge
\pi^*(K_X)|_{S_1}\cdot G\ge\frac{1}{3}(p_g(X)-4).$$
Finally, we have
$$\align
K_X^3&=\pi^*(K_X)^3\ge \pi^*(K_X)^2\cdot S_1\\
&= \pi^*(K_X)|_{S_1}\cdot S_1|_{S_1}+\pi^*(K_X)|_{S_1}\cdot E'|_{S_1}\\
&\ge (p_g(X)-2)+\frac{1}{3}(p_g(X)-4)=\frac{2}{3}(2p_g(X)-5).
\endalign$$
The inequality is sharp by virtue of (0.1). The proof is complete.
\qed
\enddemo

\remark{Remark 4.4} As was pointed out by M. Reid (\cite{R3, Remark (0.4)(v)}),
the blow-up of a canonical singularity need not be normal and thus it need not be canonical, even if the original
canonical point is a hypersurface singularity of multiplicity $2$. Because of this reason, we would rather treat a
smooth 3-fold in Theorem 4.3, although the method might be all right for Gorenstein 3-folds.
\endremark

\proclaim{Lemma 4.5} Let $X$ be a smooth projective 3-fold of general type. Suppose $p_g(X)\ge 3$,
$\dim\fei{1}(X)=1$. Keep the same notations as in subsection 1.3. If $(K_{F_0}^2, p_g(F))=(1,2),$ then one of the
following holds:

(i) $b=1$, $q(X)=1$ and $h^2(\Co{X})=0$;

(ii) $b=0$, $q(X)=0$ and $h^2(\Co{X})\le 1$.
\endproclaim
\demo{Proof}
Replacing $X$ by a birational model, if necessary, we may suppose that $\fei{1}$ is a morphism.
Note that we
do not need here the minimality of $X$. Taking the Stein-factorization of $\fei{1}$, we get a derived fibration $f:X\lrw W$.
Let $F$ be a general fiber of $f$. By assumption, $(K_{F_0}^2, p_g(F))=(1,2)$ where $F_0$ is the minimal model of
$F$. According to \cite{Ch2, Theorem 1}, we see that $b=g(W)\le 1$ whenever $p_g(X)\ge 3$. Because
$q(F)=0$, we can easily see that $q(X)=b$ and $h^2(\Co{X})=h^1(W, f_*\omega_X)$. In order to prove the lemma, It is
sufficient to study $h^1(W, f_*\omega_X)$. Since we are in a very special situation, we should be able to obtain
much more explicit information.

Let ${\Cal L}_0$ be the saturated sub-bundle of $f_*\omega_X$ which is generated by $H^0(W, f_*\omega_X)$. Because $|K_X|$ is composed of a pencil of surfaces and $\fei{1}$ factors through $f$, we see that ${\Cal L}_0$ is a line bundle on $W$. Denote ${\Cal L}_1:=f_*\omega_X/{\Cal L}_0$. Then we have the exact sequence:
$$0\lrw {\Cal L}_0\lrw f_*\omega_X\lrw {\Cal L}_1\lrw 0.$$
Noting that $\text{rk}(f_*\omega_X)=2$, we see that ${\Cal L}_1$ is also a line bundle. Noting that $H^0(W, {\Cal L}_0)\cong H^0(W, f_*\omega_X)$, we have
$h^1(W, {\Cal L}_0)\ge h^0(W, {\Cal L}_1)$. When $b=1$, \newline $\deg({\Cal L}_0)=p_g(X)\ge 3$. When $b=0$, $\deg({\Cal L}_0)=p_g(X)-1\ge 2$. Anyway, we have $h^1(W, {\Cal L}_0)=0$. So $h^0(W, {\Cal L}_1)=0$. On the other hand, It is well-known that $f_*\omega_{X/W}$ is semi-positive. Thus
$\deg({\Cal L}_1\otimes\omega_W^{-1})\ge 0$. This means $\deg({\Cal L}_1)\ge 2(b-1).$ Using the R-R, we may easily deduce that
$h^1({\Cal L}_1)\le 1-b$. So
$$h^1(W, f_*\omega_X))\le h^1(W, {\Cal L}_0)+h^1(W, {\Cal L}_1)\le 1-b.$$
So $h^2(\Co{X})\le 1-b$. The proof is complete.
\qed\enddemo

\proclaim{Lemma 4.6} Let $X$ be a smooth projective 3-fold of general type. Suppose $p_g(X)\ge 3$, $\dim\fei{1}(X)=1$ and $(K_{F_0}^2, p_g(F))=(1,2).$
Let $f:X\lrw W$ be a derived fibration of $\fei{1}$. Suppose $F_1$ and $F_2$  are two fixed smooth fibres of $f$ such that $\fei{1}(F_1)\neq \fei{1}(F_2)$.
Then
$\dim\Fi{K_X+F_1+F_2}(X)=2$  and $\Fi{K_X+F_1+F_2}|_F=\Fi{K_F}$ for a general fiber $F$.
\endproclaim
\demo{Proof}
(i). If $b=1$, we have $h^2(\Co{X})=0$ by Lemma 4.5.
       {}From the exact sequence
$$H^0(X, K_X+F_1+F_2)\lrw H^0(F_1, K_{F_1})\oplus H^0(F_2,K_{F_2})\lrw 0,$$
one may easily see that  $\dim\Fi{K_X+F_1+F_2}(X)=2$. Thus, for a general fiber $F$, $\dim\Fi{K_X+F_1+F_2}(F)=1$.
Since $p_g(F)=2$, one sees that $\Fi{K_X+F_1+F_2}|_F=\Fi{K_F}$.

(ii). If $b=0$, we only have to study $|K_X+2F_1||_F$ for a general fiber $F$. {}From the short
exact sequence:
$$0\lrw\Co{X}(K_X+F_1-F)\lrw\Co{X}(K_X+F_1)\lrw \Co{F}(K_F)\lrw 0,$$
we have the long exact sequence
$$\align
\cdots&\lrw H^0(X, K_X+F_1)\overset{\alpha_1}\to\lrw H^0(F, K_F)
\overset{\beta_1}\to\lrw H^1(X,K_X)\\
&\lrw H^1(X, K_X+F_1)\lrw H^1(F, K_F)=0,\endalign$$
If $\alpha_1$ is surjective for general $F$, then
we see that
$$\dim\Fi{K_X+F_1}(F)=\dim\Fi{K_F}(F)=1\ \text{and}\
\dim\Fi{K_X+F_1}(X)=2.$$
So $\dim\Fi{K_X+2F_1}(X)=2$. We are done. Otherwise,
$\alpha_1$ is not surjective. Because $\alpha_1\ne 0$, we see that $h^2(\Co{X})=h^1(X, K_X)\ge 1$. Because $h^2(\Co{X})\le 1$, $h^2(\Co{X})=1$ and $\beta_1$ is surjective. Therefore $H^1(X, K_X+F_1)=0$. This also means that $H^1(X, K_X+F')=0$ for any smooth fiber $F'$ since $F'\sim F_1$. So we have $H^1(X, K_X+2F_1-F)=0$, which means
$|K_X+2F_1||_F=|K_F|$. So $\dim\Fi{K_X+2F_1}(X)=2$. The proof is complete.
\qed\enddemo

\proclaim{Theorem 4.7} Let $X$ be a smooth projective 3-fold with ample canonical divisor. Suppose
$\dim\fei{1}(X)=1$ and $X$ is canonically fibered by surfaces with invariants $(c_1^2, p_g)=(1,2)$. Then
$K_X^3\ge \frac{2}{3}(2p_g(X)-7).$
\endproclaim
\demo{Proof}
The proof is slightly longer, however with the same flavour as that of Theorem 4.3.

Denote by $\overline{F}$ a generic irreducible element of $|K_X|$.
We see that $\overline{F}^2$ is a 1-cycle on $X$. If the movable part of $|K_X|$ has base points, then
$\overline{F}^2$ is a non-trivial effective 1-cycle. So $K_X\cdot \overline{F}^2\ge 2$. Thus $K_X^3\ge 2p_g(X)-2$.
Therefore we only have to treat the case when $\fei{1}$ is a morphism.

We suppose $p_g(X)\ge 3$. We still assume that $f:X\lrw W$ is a derived fibration of $\fei{1}$. Note that
$b=g(W)\le 1$.  Let $\overline{M}$ be the movable part of $|K_X+F_1+F_2|$.
 Also note that $F$ is minimal in this situation and $(K_F^2, p_g(F))=(1,2)$. It is well-known that $|K_F|$
has exactly one base point, but no fixed part, and that a general member of $|K_F|$ is a smooth irreducible curve of
genus 2. Since $|K_X+F_1+F_2||_F=|K_F|$ and according to Lemma 2.6, we see that $\overline{M}|_F=K_F$. This
means that $|\overline{M}|$ definitely has base points.
According to Hironaka, we can take successive blow-ups
$$\pi: X'=X_{n}\overset{\pi_n}\to\rw X_{n-1}\rw\cdots\rw
   X_{i}\overset{\pi_i}\to\rw X_{i-1}\rw\cdots\rw
  X_1\overset{\pi_1}\to \rw X_0=X$$
such that

(i) $\pi_i$ is a single blow-up along smooth center $W_i$ on $X_{i-1}$ for all $i$;

(ii) $W_i$ is contained in the base locus of the movable part of
$$|(\pi_1\circ\pi_2\circ\cdots\circ\pi_{i-1})^*(\overline{M})|$$
and thus $W_i$ is a reduced closed point or a smooth projective curve on $X_{i-1}$;

(iii) the movable part of $|\pi^*(\overline{M})|$ has no base points.

Denote by $E_i$ the exceptional divisor on $X'$ corresponding to $W_i$ for all $i$. Note that the resulting 3-fold $X'$ is still smooth. Let $M$ be the movable part of $|\pi^*(\overline{M})|$ and $S\in|M|$ be a general member. Then $S$ is a smooth irreducible projective surface of general type. Denote $f':=f\circ\pi$. Then $f':X'\lrw W$ is still a fibration. Let $F'$ be a general fiber of $f'$. Note that $F'$ has the minimal model $F$. We may write
$$K_{X'}\sim \pi^*(K_X)+\sum_{i=0}^pa_iE_i=\pi^*(K_X)+E$$
and $\pi^*(\overline{M})=M+\sum_{i=0}^pe_iE_i$. According to Lemma 4.2, we have $0<a_i\le 2e_i$ for all $i$.
Recall that we have
$K_X\sim S_1+Z=\sum_{i=1}^{b_1}F_i+Z,$
where $b_1\ge p_g(X)-1$, the $F_i's$ are fibers of $f$, $S_1$ is the movable part of $|K_X|$ and $Z$ the fixed part of $|K_X|$. Note that there is an effective divisor $Z_0\le Z$ such that
$\overline{M}\sim S_1+F_1+F_2+Z_0.$
We  write
$$\align
\pi^*(K_X+F_1+F_2)&\sim \pi^*(\overline{M}+Z-Z_0)=M+\sum_{i=0}^pe_iE_i+
\pi^*(Z-Z_0)\\
&=M+\sum_{i=0}^pe_i'E_i+\sum_{j=1}^qd_jL_j=:M+E',
\endalign$$
where $E_i\ne L_j$, $d_j>0$, $e_i'\ge e_i$ for all $i$ and $L_{j_1}\ne L_{j_2}$ whenever $j_1\ne j_2$.
Note that $\pi^*(\overline{M})\ge \pi^*(S_1+F_1+F_2)$ and that the strict transform of $S_1$ is a union of $b_1$ fibers of $f'$, we see that
$$M|_S\ge \sum_{j=1}^{b_1+m}F_j'|_S\equiv (b_1+m)F'|_S$$
where the $F_j's$ are fibers of $f'$ and $m=2$.
Because $\dim\Fi{M}(X')=2$, we see $\dim\Fi{M}(S)=1$ for a general member $S$.
So, on $S$, the system $|M|_S|$ should be composed of a free pencil of curves since $(M|_S)^2=M^3=0$. On the other hand, we obviously have
$H^0(X', K_{X'}-S)=0$. This instantly gives the inclusion
$H^0(X', K_{X'})\hookrightarrow H^0(S, K_{X'}|_S).$ So $\dim\Fi{K_{X'}}(S)\ge 1$.
Because $\dim\fei{1}(X)=1$,  we see that $\dim\Fi{K_{X'}}(S)=1$. Thus It is clear
$f'(S)=W$. So we have a surjective morphism $f'|_S: S\lrw W$. The fiber of $f'|_S$ is exactly $F'\cap S$ or the divisor $F'|_S$. Since $|M|_S|$ is composed of a pencil of curves,
$M|_S\ge\sum_{j=1}^{b_1+m}F_j'|_S$ and $|\sum_{j=1}^{b_1+m}F_j'|_S|$ is vertical, we see that $|M|_S|$ is also
vertical, i.e.
$\dim f'|_S(M|_S)=0$. This means that the divisor $M|_S$ is vertical with respect to the morphism $f'|_S$. By taking the Stein-factorization of $f'|_S$, one can see that $F'|_S$ is linearly equivalent to a disjoint union of irreducible curves of the same numerical type and  $F'|_S\equiv\xi C$ where $C$ is certain irreducible curve and $\xi$ is a positive integer.

Recall that $E':=\sum_{i=0}^pe_i'E_i+\sum_{j=1}^qd_jL_j$. We may write
$E'|_S:=E_V'+E_H'$ where $E_V'$ is the vertical part and $E_H'$ is the horizontal part with $E_H'\cdot F'|_S>0$. Noting that $\pi^*(K_X+F_1+F_2)|_S$ is nef and big and that $M|_S$ is vertical, we see that $E_H'$ is non-trivial. So we have
$$\pi^*(K_X+F_1+F_2)|_S=M|_S+E'|_S=M|_S+E_V'+E_H'.$$
Also recall that $E:=\sum_{i=0}^pa_iE_i$. Denote $E|_S:=E_V+E_H$ where $E_V$ is the vertical part and $E_H$
is the horizontal part. We have
$$\align
0&<F'|_S\cdot E_H'=F'|_S\cdot E'|_S=F'|_S\cdot \pi^*(K_X+F_1+F_2)|_S\\
&=F'\cdot \pi^*(K_X+F_1+F_2)\cdot S\\
&\le F'\cdot \pi^*(K_X+F_1+F_2)\cdot \pi^*(K_X+F_1+F_2)=
K_X^2\cdot F=1.\endalign$$
This means
$$F'|_S\cdot E_H'=F'|_S\cdot \pi^*(K_X)|_S=1, \tag{4.5}$$
$$\pi^*(F_1)|_S\cdot F'|_S=0. \tag{4.6}$$
Thus we see that $\xi=1$ and thus $f'|_S: S\lrw W$ is a fibration.
This also means that $E_H'$ is irreducible and that it comes from certain irreducible component of $E'$. For generic
$S$ and $F'$, because $S|_{F'}$ is the movable part of $|K_{F'}|$, we see that $S|_{F'}$ is an irreducible curve of genus two. This means $C=S\cap F'$ is a smooth curve of genus
2 on $S$ and $C^2=(F'|_S)^2=0$. Thus $K_S\cdot C=2$, i.e.
$$(E_V+E_H+\pi^*(K_X)|_S+S|_S)\cdot C=2.$$
Noting that, from (4.5), $S|_S\cdot C=M|_S\cdot F'|_S=0$ and $\pi^*(K_X)\cdot C=1$, we have $E_H\cdot C=1$. This also says that $E_H$ comes from certain irreducible component $E_i$ in $E$ with $a_i=1$. For simplicity we may suppose that this component is just $E_0$. So $a_0=1$. Now It is quite clear about the structure of $E'|_S$ and $E|_S$:
$$E_H=E_H'\le E_0|_S,\ \ \sum_{i=1}^p a_i(E_i|_S)+(E_0|_S-E_H)=E_V,$$
$$\sum_{i=1}^pe_i'(E_i|_S)+\sum_{j=1}^qd_j(L_j|_S)+(E_0|_S-E_H')=E_V'.$$
Noting that $E_0|_S$ can have only one horizontal component, we denote it by
$G:=E_H=E_H'$. Similar to the Claim in the proof of Theorem 4.3, It is easy to see that $E_V\le 2E_V'$.

Now we may perform the computation on the surface $S$. We have
$$(K_S+G+2(1-b)F'|_S)\cdot G=2p_a(G)-2+2(1-b)\ge 0.$$
(One notes that $p_a(G)\ge 1$ if $b=1$ and $p_a(G)\ge 0$ if $b=0$.)
$$\align
K_S\cdot G&=(E|_S+\pi^*(K_X)|_S+S|_S)\cdot G=
E_V\cdot G+G^2+\pi^*(K_X)|_S\cdot G+S|_S\cdot G\\
&\le 2E_V'\cdot G+G^2+S|_S\cdot G+\pi^*(K_X)|_S\cdot G\\
&=E_V'\cdot G+\pi^*(K_X+F_1+F_2)|_S\cdot G+\pi^*(K_X)|_S\cdot G.
\endalign$$
So we get
$$E_V'\cdot G+\pi^*(2K_X+F_1+F_2)|_S\cdot G +G^2+2(1-b)\ge 0. \tag{4.7}$$
On the other hand, we have
$$\align
\pi^*(K_X+F_1+F_2)|_S\cdot G&=S|_S\cdot G+E_V'\cdot G+G^2\\
&\ge (b_1+m)F'|_S\cdot G+E_V'\cdot G+G^2, \tag{4.8}
\endalign$$
where we note that $S|_S$ is vertical and, numerically,  $S|_S\ge_{\text{num}} (b_1+m)F'|_S$
and $F'|_S\cdot G=1$ by (4.5).
Combining (4.7) and (4.8), we get
$$\pi^*(3K_X+2F_1+2F_2)|_S\cdot G\ge (b_1+m)+2(b-1).$$
We have
$$\align
&\pi^*(3K_X+2F_1+2F_2)|_S\cdot G \le \pi^*(3K_X+2F_1+2F_2)|_S\cdot E'|_S\\
=&\pi^*(3K_X+2F_1+2F_2)|_S\cdot (\pi^*(K_X+F_1+F_2)|_S-S|_S)\\
=&\pi^*(3K_X+2F_1+2F_2)|_S\cdot \pi^*(K_X+F_1+F_2)|_S-\pi^*(3K_X+2F_1+2F_2)|_S\cdot S|_S\\
\le &(3K_X+2F_1+2F_2)(K_X+F_1+F_2)^2-\pi^*(3K_X+2F_1+2F_2)|_S\cdot S|_S\\
= &3K_X^3+8m-\pi^*(3K_X+2F_1+2F_2)|_S\cdot S|_S.
\endalign$$
Thus
$3K_X^3\ge b_1-7m+2(b-1)+\pi^*(3K_X+2F_1+2F_2)|_S\cdot S|_S.$
By (4.5) and (4.6), we get
$$\pi^*(3K_X+2F_1+2F_2)|_S\cdot S|_S\ge \pi^*(3K_X+2F_1+2F_2)|_S\cdot (b_1+m)F'|_S=
3(b_1+m).$$
So $3K_X^3\ge 4b_1-4m+2(b-1).$
We obtain
$$K_X^3\ge \frac{4}{3}b_1-\frac{4}{3}m+\frac{2}{3}(b-1)\ge
\cases \frac{4}{3}p_g(X)-\frac{8}{3}, \ &\ \text{if}\  b=1\\
\frac{4}{3}p_g(X)-\frac{14}{3}, \  &\ \text{if}\  b=0.\endcases$$

Finally, we discuss what happens when $K_X^3>\frac{4}{3}p_g(X)-\frac{10}{3}$.
Definitely, $b=0$ and $3K_X^3=4p_g(X)-11$, $4p_g(X)-12$, $4p_g(X)-13$ or $4p_g(X)-14$.
Noting that $K_X^3$ is an even number, one excludes possibilities $4p_g(X)-11$ and $4p_g(X)-13$.
The proof is complete.
\qed\enddemo

\proclaim{Corollary 4.8} Let $X$ be a smooth projective 3-fold with ample canonical divisor. Then
we have the following  Noether inequality
$$K_X^3\ge\frac{2}{3}(2p_g(X)-7).$$
\endproclaim
\demo{Proof}
This is a direct result of Theorem 4.1, Theorem 4.3 and Theorem 4.7.
\qed\enddemo

Corollary 4.8 implies Corollary 2. Theorem 4.1, Theorem 4.3 and Theorem 4.7 imply Theorem 5(1) and
Theorem 5(2).

\head {\rm 5. An appendix}\endhead
We go on proving Theorem 5 in this section.

\proclaim{Proposition 5.1}
Let $X$ be a projective minimal Gorenstein 3-fold of general type with only locally factorial terminal singularities.
Suppose $X$ has a locally factorial canonical model.
If $\dim\fei{1}(X)=1$and  $(K_{F_0}^2, p_g(F))=(1,2)$,  then
$$K_X^3\ge \frac{2}{21}(11p_g(X)-16).$$
\endproclaim
\demo{Proof}
If the movable part of $|K_X|$ has base points, then we have $K_X^3\ge 2p_g(X)-2$ according to
\cite{Kob, Case 1, Theorem (4.1)} because $X$ is assumed to have a locally factorial canonical model.
So we may suppose $\Fi{K_X}$ is a morphism.

Taking the Stein-factorization of $\Fi{K_X}$, we get the derived fibration
$f:X\lrw W$. Let $M_1$ be the movable part of $|K_X|$ and $S_1\in|M_1|$ a general member. We may write
$S_1\sim\sum_{i=1}^{b_1}F_i\equiv b_1F,$
where the $F_i's$ are fibers of $f$, $F$ is a general fiber of $f$ and $b_1\ge p_g(X)-1$. Because $X$ is minimal, $F$ is a minimal surface.
Since $X$ has isolated singularities, $F$ is smooth.
Note that we have $K_F^2=1$ and $p_g(F)=2$  under the assumption of the proposition. We may also write
$K_X\equiv b_1F+Z$, where $Z$ is the fixed part of $|K_X|$. According to \cite{Ch2, Theorem 1}, we have
$b:=g(W)\le 1$ provided $p_g(X)\ge 3$. {}From \cite{L}, we know that $|4K_X|$ is base point free.
Let $S_4\in |4K_X|$ be a general member. Since $X$ has isolated singularities, $S_4$ is a smooth projective
irreducible surface of general type.  we see that $f(S_4)=W$. Denote $f_0:=f|_{S_4}$. Then $f_0: S_4\lrw W$ is a proper surjective
morphism onto $W$ ($f_0$ need not be a fibration). Because
$f(F)$ is a point, $F|_{S_4}$ is vertical with respect to $f_0$, {\it i.e.,}
$\dim f_0(F|_{S_4})=0$. Now we have
$K_X|_{S_4}\equiv b_1F|_{S_4}+Z|_{S_4}.$
Denote $Z|_{S_4}:=Z_V+Z_H$, where $Z_V$ is the vertical part and $Z_H$ is the horizontal part. We may write
$Z_H:=\sum m_iG_i$, where $m_i>0$ and the $G_i's$ are distinct irreducible curves on $S_4$. We have
$$\align
&(F|_{S_4}\cdot Z_H)_{S_4}=(F|_{S_4}\cdot Z|_{S_4})_{S_4}=(F\cdot S_4\cdot Z)_X\\
=&(S_4|_F\cdot Z|_F)_F=4(K_X|_F\cdot K_X|_F)_F=4K_F^2=4.
\endalign$$
Thus $m_i\le 4$ for all $i$. Denote
$$D:=4K_{S_4}-8(b-1)F|_{S_4}+Z_V+Z_H.$$
We claim that $D\cdot G_i\ge 0$ for all $i$. In fact, since $Z_V\cdot G_i\ge 0$
and $G_i\cdot G_j\ge 0$ for $i\ne j$, we have
$$\align
D\cdot G_i&\ge 4K_{S_4}\cdot G_i-8(b-1)F|_{S_4}\cdot G_i+m_iG_i^2\\
&=(4-m_i)K_{S_4}\cdot G_i+m_i(K_{S_4}\cdot G_i+G_i^2)-8(b-1)F|_{S_4}\cdot G_i\\
&=(4-m_i)K_{S_4}\cdot G_i+m_i(2p_a(G_i)-2)-8(b-1)F|_{S_4}\cdot G_i.
\endalign$$
Note that both $K_{S_4}$ and $F|_{S_4}$ are nef. When $b=1$, we have $p_a(G_i)\ge b=1$.
Thus $D\cdot G_i\ge (4-m_i)K_{S_4}\cdot G_i\ge 0$. When $b=0$,
$$D\cdot G_i\ge (4-m_i)K_{S_4}\cdot G_i+(8-2m_i)F|_{S_4}\cdot G_i+m_i[2p_a(G_i)-2+2F|_{S_4}\cdot G_i]\ge 0.$$
Therefore we get $D\cdot Z_H\ge 0$. This means
$$4K_{S_4}\cdot Z_H-8(b-1)F|_{S_4}\cdot Z_H+(Z_V+Z_H)Z_H\ge 0. \tag{5.1}$$
On the other hand, we have
$$K_X|_{S_4}\cdot Z_H=b_1F|_{S_4}\cdot Z_H+(Z_V+Z_H)Z_H. \tag{5.2}$$
Combining (5.1) and (5.2), we get
$$\align
4K_{S_4}\cdot Z_H+K_X|_{S_4}\cdot Z_H&\ge (b_1+8(b-1))F|_{S_4}\cdot Z_H\\
&\ge 4(p_g(X)+10b-11).
\endalign$$
We also have
$$\align
&4K_{S_4}\cdot Z_H+K_X|_{S_4}\cdot Z_H=5K_X|_{S_4}\cdot Z_H+4S_4|_{S_4}\cdot Z_H\\
\le&5K_X|_{S_4}\cdot Z|_{S_4}+4S_4|_{S_4}\cdot Z|_{S_4}=84K_X^2\cdot Z.
\endalign$$
Thus we obtain
$$K_X^2\cdot Z\ge \frac{1}{21}(p_g(X)+10b-11)=\cases
\frac{1}{21}(p_g(X)-11),\ &\text{if}\ b=0,\\
\frac{1}{21}(p_g(X)-1),\ &\text{if}\ b=1.\endcases $$
Finally we get
$$
K_X^3\ge b_1K_X^2\cdot F+K_X^2\cdot Z
\ge\cases
\frac{2}{21}(11p_g(X)-16), \ &\text{if}\ b=0,\\
\frac{22}{21}(p_g(X)-1), \ &\text{if}\ b=1.\endcases$$
The proof is complete.
\qed
\enddemo

Section 4 and Proposition 5.1 imply Theorem 5(3).

\head  {\rm Acknowledgment }\endhead
This paper was begun when the author was a postdoctor at Georg-August Universitat Goettingen,
Germany between October 1999 and August 2000.  The author thanks F. Catanese and M. Reid for
generous helps. Especially this note benefited a lot from frequent discussions between
professor Catanese and the author  while at Goettingen.
 Special thanks
are due to both M. Kobayashi who kindly answered the author's query on his paper through email
and Seunghun Lee who pointed out an error to the author in the proof of Lemma 4.2 in the first
draft of the paper. Finally the author heartly appreciates the referee's invaluable suggestions
which essentially make this note more readable and better organized.

\head {\rm References} \endhead
\roster
\item"[Be1]" A. Beauville, {\it L'application canonique pour les surfaces de type g\'en\'eral}, Invent. Math. {\bf 55}(1979),
121-140.
\item"[Bo]" E. Bombieri, {\it Canonical models of surfaces of general type},
Publications I.H.E.S. {\bf 42}(1973), 171-219.
\item"[Cas]" G. Castelnuovo, {\it Osservazioni intorno alla geometria sopra una superficie, I, II}, Rendiconti del R. Istituto Lombardo, s. II, {\bf 24}(1891), also in `Memorie scelte', Zanichelli (1937), Bologna, 245-265.
\item"[Cat1]" F. Catanese, {\it Surfaces with $K^2=p_g=1$ and their period mapping}, Springer Lecture Notes in
Math. {\bf 732}(1979), pp. 1-29.
\item"[Ch1]" M. Chen, {\it An inequality with regard to the canonical birationality for varieties}, Commun. Algebra {\bf 23}(1995), 4439-4446.
\item"[Ch2]" ------, {\it Complex varieties of general type whose canonical systems are composed with pencils}, J. Math. Soc. Japan {\bf 51}(1999), 331-335.
\item"[Ch3]" ------, {\it Kawamata-Viehweg vanishing and quint-canonical maps for threefolds of general type}, Commun. Algebra {\bf 27}(1999), 5471-5486.
\item"[Ch4]" ------, {\it Canonical stability in terms of singularity index for algebraic threefolds}, Math. Proc. Camb. Phil. Soc. {\bf 131}(2001), 241-264.
\item"[Ci]" C. Ciliberto, {\it The bicanonical map for surfaces of general type}, Proc. Symposia in Pure Math. {\bf 62}(1997), 57-83.
\item"[E-L]" L. Ein, R. Lazarsfeld, {\it Global generation of pluricanonical
and adjoint linear systems on smooth projective threefolds}, J. Amer. Math. Soc.
{\bf 6}(1993), 875-903.
\item"[F]" T. Fujita, {\it On Kahler fiber spaces over curves}, J. Math. Soc. Japan {\bf 30} (1978), 779-794.
\item"[Hrs]" J. Harris, {\it A bound on the geometric genus of projective varieties}, Ann. Sc. Norm. Sup. Pisa {\bf 8}(1981), 35-68.
\item"[Hi]" H. Hironaka, {\it Resolution of singularities of an algebraic variety over a field of characteristic zero}, I, Ann. of Math. {\bf 79}(1964), 109-203, II, ibid., 205-326.
\item"[Ka1]" Y. Kawamata, {\it A generalization of Kodaira-Ramanujam's
vanishing theorem}, Math. Ann. {\bf 261}(1982), 43-46.
\item"[KMM]" Y. Kawamata, K. Matsuda, K. Matsuki, {\it Introduction to the
minimal model problem}, Adv. Stud. Pure Math. {\bf 10}(1987), 283-360.
\item"[Kob]" M. Kobayashi, {\it On Noether's inequality for threefolds}, J. Math. Soc. Japan {\bf 44}(1992), 145-156.
\item"[Kol]" J. Koll\'ar, {\it Higher direct images of dualizing sheaves}, I,
Ann. of Math. {\bf 123}(1986), 11-42;\ \ II, ibid. {\bf 124}(1986), 171-202.
\item"[L]" S. Lee,  {\it Quartic-canonical systems on canonical threefolds of index 1},
Commun. Algebra   {\bf 28}(2000), 5517-5530.
\item"[M1]" Y. Miyaoka, {\it On the chern numbers of surfaces of general type}, Invent. Math. {\bf 42}(1977), 225-237.
\item"[M2]" ------, {\it The pseudo-effectivity of $3c_2-c_1^2$ for varieties with numerically effective canonical classes}, Algebraic Geometry, Sendai, 1985.
Adv. Stud. Pure Math. {\bf 10}(1987), 449-476.
\item"[M3]" ------, {\it On the Mumford-Ramanujam vanishing theorem on a surface}.
\newline Journ\'ees de G\'eom\'etrie
Alg\'ebrique d'Angers (A. Beauville, ed.), Sijthoff and Noordhof, Alphen ann den
Rijn 1980, 239-247.
\item"[N]" M. Noether, {\it Zur Theorie des eindeutigen Entsprechens algebraischer Ge-bilde}, Math. Ann. {2}(1870), 293-316; {\bf 8}(1875), 495-533.
\item"[R1]" M. Reid, {\it Young person's guide to canonical singularities}, Proc. Symposia in pure Math. {\bf 46}(1987), 345-414.
\item"[R2]" ------, {\it Minimal models of canonical 3-folds}, Algebraic Varieties and Analytic Varieties (S. Itaka, ed.), Adv. Stud. Pure Math. {\bf 1}(1983), 131-180.
\item"[R3]" ------, {\it Canonical 3-folds}, Journ\'ees de G\'eom\'etrie
Alg\'ebrique d'Angers (A. Beauville, ed.), Sijthoff and Noordhof, Alphen ann den
Rijn 1980, 273-310.
\item"[Rr]" I. Reider, Vector bundles of rank 2 and linear systems on
algebraic surfaces, {\it Ann. Math.} {\bf 127}(1988), 309-316.
\item"[S]" F. Sakai, {\it Weil divisors on normal surfaces}, Duke Math. J.
 {\bf 51}(1984), 877-887.
\item"[V1]" E. Viehweg, {\it Vanishing theorems}, J. reine angew. Math. {\bf
335}(1982), 1-8.
\item"[V2]" ------, {\it Weak positivity and the additivity of the Kodaira dimension for certain fibre spaces}, Adv. Stud. Pure Math. {\bf 1}(1983), 329-353.
\item"[V3]" ------, {\it Weak positivity and the additivity of the Kodaira dimension, II: The local Torelli map}, Classification of algebraic and analytic manifolds, Prog. Math. {\bf 39}(1983), 567-589.
\item"[Y1]" S. T. Yau, {\it Calabi's conjecture and some new results in algebraic geometry}, Proc. Nat. Acad. Sci. USA {\bf 74}(1977), 1798-1799.
\item"[Y2]" ------, {\it On the Ricci curvature of a complex K$\ddot{\text{a}}$hler manifold and the complex Monge-Ampere equations}, Commun.
Pure and Appl. Math. {\bf 31}(1978), 339-411.
\endroster
\enddocument